\documentclass[12pt]{article}
\usepackage{amssymb}
\usepackage{amsmath}
\numberwithin{equation}{section}

\newcommand{\mb}[1]{\ \mbox{ #1 }\ }
\newcommand{\nonu}{\nonumber \\}

\newtheorem{prop}{Proposition}[section]

\newtheorem{defi}[prop]{Definition}
\newtheorem{theo}[prop]{Theorem}

\newcommand{\prf}{\underline{Proof:}\ }
\newcommand{\finprf}{\null \hfill {\rule{5pt}{5pt}}\\[2.1ex]\indent}

\newcommand{\RR}{\mbox{${\mathbb R}$}}

\newcommand{\ZZ}{\mbox{${\mathbb Z}$}}

\newcommand{\II}{\mbox{${\mathbb I}$}}
\newcommand{\ca}{\mbox{$\cal{A}$}}
\newcommand{\cb}{\mbox{$\cal{B}$}}
\newcommand{\cc}{\mbox{$\cal{C}$}}

\newcommand{{\cg}}{\mbox{$\cal{G}$}}

\newcommand{\ci}{\mbox{${\cal I}$}}

\newcommand{\ck}{\mbox{$\cal{K}$}}

\newcommand{\cR}{\mbox{$\cal{R}$}}

\newcommand{\ct}{\mbox{$\cal{T}$}}
\newcommand{\cu}{\mbox{${\cal U}$}}

\newcommand{\invdots}{{\mathinner{\mkern1mu\raise1pt
     \vbox{\kern7pt\hbox{.}}\mkern2mu\raise4pt\hbox{.}
     \mkern2mu\raise7pt\hbox{.}\mkern1mu}}}

  \newcommand{\1}{\mathbf{1}}
\newcommand{\half}{\frac{1}{2}}
\newcommand{\ie}{{\it i.e.}\ }
\begin{document}
\renewcommand{\thefootnote}{\fnsymbol{footnote}}
\newpage
\pagestyle{empty}
\setcounter{page}{0}
\pagestyle{empty}
\newcommand{\LAP}{LAPTH}
\def\logo{{\bf {\huge LAPTH}}}

\centerline{\logo}

\vspace {.3cm}

\centerline{{\bf{\it\Large 
Laboratoire d'Annecy-le-Vieux de Physique Th\'eorique}}}

\centerline{\rule{12cm}{.42mm}}

\vspace{20mm}

\begin{center}

  {\LARGE  {\sffamily Interplay between Zamolodchikov-Faddeev
  \\[1.2ex]
 and Reflection-Transmission algebras}}
\\[2.1em]

{\large
M. Mintchev$^{a}$\footnote{mintchev@df.unipi.it} and
E. Ragoucy$^{b}$\footnote{ragoucy@lapp.in2p3.fr}
}\\
\end{center}

\null

\noindent
{\it $^a$ INFN and Dipartimento di Fisica, Universit\'a di
     Pisa, Via Buonarroti 2, 56127 Pisa, Italy\\[2.1ex]
$^b$ LAPTH, 9, Chemin de Bellevue, BP 110, F-74941 Annecy-le-Vieux
     cedex, France}
\vfill

\begin{abstract}
We show that a suitable coset algebra,
constructed in terms of an extension of the Zamolodchikov-Faddeev algebra, is homomorphic to the 
Reflection-Transmission algebra, as it appears in the study of 
integrable systems with impurity.
\end{abstract}

\vfill
\centerline{J. Phys. \textbf{A37} (2004) 425}
\vfill
\rightline{May 2003}
\rightline{IFUP-TH 18/2003}
\rightline{LAPTH-983/03}
\rightline{\tt math.QA/0306084}
\newpage

\markright{\today\dotfill DRAFT\dotfill }
\pagestyle{myheadings}
\pagestyle{plain}

\setcounter{footnote}{0}

\section{Introduction}

The Zamolodchikov-Faddeev (ZF) algebra \cite{Zamolodchikov:xm, Faddeev:zy}
is well-known to be the basis
of factorized scattering theory in  integrable models on the
line. The two-body scattering, which is the only particle interaction
present in this case, is implemented by a quadratic constraint among
the particle creation and annihilation operators. There exist various
attempts \cite{Cherednik:vs}-\cite{Castro-Alvaredo:2002dj}
to generalize factorized scattering theory
to the  case when point-like impurities, preserving integrability, are present
as well. In this context, the relevant algebraic structure emerging recently
\cite{Mintchev:2002zd, Mintchev:2003ue}, is the so-called  reflection-transmission (RT) algebra. Besides the two-body scattering
in the bulk, the RT algebra captures also the particle interaction 
with the impurity. From
the physical analysis performed in \cite{Mintchev:2002zd},
it follows that the particle-impurity
interactions are implemented  by further constraints, ensuring the
compatibility  between bulk scattering and the
process of reflection/transmission from the impurity.
This feature suggest the existence of an algebraic connection 
among ZF and RT algebras.
It turns out in fact that taking the coset of the ZF algebra with respect to
a suitable two-sided ideal, one obtains an algebra homomorphic
to the RT algebra. The explicit realization of this relationship
is the main goal of this investigation. Our analysis
extends the results of \cite{Ragoucy:2001cf, ZFconf}, obtained in the 
case of pure reflection.
Establishing the link between ZF and RT algebras, we take 
the occasion to discuss
also the issue of symmetries when impurities are present.

\section{Background}

We collect here some definitions and basic results concerning ZF and 
RT algebras,
which are needed in what follows. Adopting throughout the paper the
compact tensor notation, introduced
in \cite{Ragoucy:2001zy}, we start by:

\begin{defi}{\bf (ZF algebra $\ca_{S}$)}\\
$\ca_{S}$ is the polynomial algebra generated by a unit 
element $\1$ and the 
generators $a(k)$
and $a^\dag(k)$, subject to the constraints:
\begin{eqnarray}
a_{1}(k_1)a_{2}^\dag(k_2) &=&
a_{2}^\dag(k_2)S_{12}(k_{1},k_{2})a_{1}(k_1)+
\delta_{12}\delta(k_1-k_2) \1 \quad\\
a_{1}(k_1)a_{2}(k_2) &=& S_{21}(k_{2},k_{1})
a_{2}(k_2)a_{1}(k_1)\\
a^\dag_{1}(k_1)a_{2}^\dag(k_2) &=&
a_{2}^\dag(k_2)a^\dag_{1}(k_1) S_{21}(k_{2},k_{1})
\end{eqnarray}
\end{defi}
Here and below the $S$-matrix obeys the following well-known
relations (Yang-Baxter equation and unitarity):
\begin{eqnarray}
&&S_{12}(k_{1},k_{2})\,S_{13}(k_{1},k_{3})\,S_{23}(k_{2},k_{3}) =
     S_{23}(k_{2},k_{3})\,S_{13}(k_{1},k_{3})\,S_{12}(k_{1},k_{2})\\
&&S_{12}(k_{1},k_{2})\,S_{21}(k_{2},k_{1}) = \II\otimes\II
\end{eqnarray}
We emphasize that $S$ depends in general on $\chi_1$ and $\chi_2$, 
observing in addition that
bulk $S$-matrices depending on $\chi_1-\chi_2$ only are also  covered 
by our framework, which therefore allows
to treat both systems with exact as well as broken 
Lorentz (Galilean) invariance. More details about the physical meaning 
of this generalization are given in \cite{Mintchev:2003ue}.

We shall employ below the extension ${\overline\ca}_{S}$ of $\ca_{S}$,
involving power series in $a(k)$ and $a^\dag(k)$ whose individual terms
preserve the particle number. The concept of extended ZF (EZF) algebra
${\overline\ca}_{S}$ is relevant for proving \cite{Ragoucy:2001zy} 
the following:
\begin{prop}{\bf (Well-bred operators)}\\
Each ${\overline\ca}_{S}$ contains an invertible element $L(k)$ satisfying
\begin{eqnarray}
L_{1}(k_{1})a_{2}(k_{2})&=&
S_{21}(k_{2},k_{1})a_{2}(k_{2})L_{1}(k_{1}) \\
L_{1}(k_{1})a^\dag_{2}(k_{2}) &=&
a^\dag_{2}(k_{2})S_{12}(k_{1},k_{2})L_{1}(k_{1})
\end{eqnarray}
Moreover, $L(k)$ obeys
\begin{equation}
S_{12}(k_{1},k_{2})L_{1}(k_{1})L_{2}(k_{2})=
L_{2}(k_{2})L_{1}(k_{1})S_{12}(k_{1},k_{2})
\end{equation}
and generates a quantum group $\cu_{S}\subset{\overline\ca}_{S}$,
with Hopf structure $\Delta L(k)= L(k)\otimes L(k)$.
$L(k)$ also satisfies $L(k)^\dag= L(k)^{-1}$.
\end{prop}
$L(k)$, called well-bred operator, is explicitly constructed in
\cite{Ragoucy:2001zy} and admits a series representation in terms of 
$a(k)$ and $a^\dag(k)$.

We turn now to RT algebras \cite{Mintchev:2003ue}.
\begin{defi}{\bf (RT algebra $\cc_{S}$)}\\
A RT algebra is generated by $\1$ and the generators $A(k)$, 
$A^\dag(k)$, $t(k)$ and $r(k)$
obeying:
\begin{eqnarray}
A_{1}(k_1)A_{2}(k_2) &=& S_{21}(k_{2},k_{1})
A_{2}(k_2)A_{1}(k_1)\label{rt-1}\\
A^\dag_{1}(k_1)A_{2}^\dag(k_2) &=&
A_{2}^\dag(k_2)A^\dag_{1}(k_1) S_{21}(k_{2},k_{1})
\label{rt-2}\\
A_{1}(k_1)A_{2}^\dag(k_2) &=&
A_{2}^\dag(k_2)S_{12}(k_{1},k_{2})A_{1}(k_1)+
\delta(k_1-k_2)\delta_{12} \Big(\1 +t_{1}(k_{1})\Big)\nonu
&& +\delta_{12} r_{1}(k_1)\delta(k_1+k_2)\quad
\label{rt-3}
\end{eqnarray}
\begin{eqnarray}
A_{1}(k_1)t_{2}(k_2) &=& S_{21}(k_{2},k_{1})
t_{2}(k_2)S_{12}(k_{1},k_{2}) A_{1}(k_1)
\label{rt-4}\\
A_{1}(k_1)r_{2}(k_2) &=& S_{21}(k_{2},k_{1})
r_{2}(k_2)S_{12}(k_{1},-k_{2}) A_{1}(k_1)
\label{rt-5}\\
t_{1}(k_1)A_{2}^\dag (k_2) &=&
A_{2}^\dag(k_2) S_{12}(k_{1},k_{2}) t_{1}(k_1)S_{21}(k_{2},k_{1})
\label{rt-6}\\
r_{1}(k_1)A_{2}^\dag (k_2) &=&
A_{2}^\dag(k_2) S_{12}(k_{1},k_{2}) r_{1}(k_1)S_{21}(k_{2},-k_{1})
\qquad
\label{rt-7}
\end{eqnarray}
\begin{eqnarray}
\!\! && \!\!\!\! S_{12}(k_{1},k_{2})\, t_{1}(k_1)\,
S_{21}(k_{2},k_{1})\, t_{2}(k_2)\, =\,
t_{2}(k_2)\, S_{12}(k_{1},k_{2})\, t_{1}(k_1)\, S_{21}(k_{2},k_{1})
\label{rt-8}\\
&& \!\!\!\! S_{12}(k_{1},k_{2})\, t_{1}(k_1)\, S_{21}(k_{2},k_{1})\,
r_{2}(k_2)\, =\,
r_{2}(k_2)\, S_{12}(k_{1},-k_{2})\, t_{1}(k_1)\, S_{21}(-k_{2},k_{1})
\label{rt-9}\\
&& \!\!\!\! S_{12}(k_{1},k_{2})\, r_{1}(k_1)\, S_{21}(k_{2},-k_{1})\,
r_{2}(k_2)\, =\,
r_{2}(k_2)\, S_{12}(k_{1},-k_{2})\, r_{1}(k_1)\, S_{21}(-k_{2},-k_{1})
\qquad\quad\label{rt-10}
\end{eqnarray}
where $r(k)$ and $t(k)$ satisfy:
\begin{eqnarray}
t(k)t(k)+r(k)r(-k) &=& \1\label{rt-12a}\\
   \qquad t(k)r(k)+r(k)t(-k) &=& 0 \label{rt-12b}
\end{eqnarray}
\end{defi}
We refer to $r(k)$ and $t(k)$ as reflection and transmission generators.

A special case of RT algebra is the boundary algebra $\cb_{S}$ 
\cite{Liguori:1998xr}
defined as the coset of $\cc_{S}$ by the relation $t(k)=0$. More
precisely, one defines the left ideal $\ci_{S}=\cc_{S}.\{t(k)\}$. From
the above relations, this left ideal is obviously a two-sided ideal,
so that the coset $\cc_{S}/\ci_{S}$ is an algebra: it is the
boundary algebra $\cb_{S}$, whose defining relations are given by
eqs. (\ref{rt-1})--(\ref{rt-12a}) with $t(k)=0$.

Note that the same construction applies when
$r(k)=0$, leading to an algebra associated to a purely transmitting
impurity \cite{Bowcock:2003dr}.

Note also that the generators $r(k)$ and $t(k)$ form a
subalgebra $\ck_{S}\subset \cc_{S}$. This subalgebra contains itself 
two subalgebras.
The subalgebra $\cR_{S}\subset \ck_{S}$, generated by $r(k)$, appears 
already in
\cite{Sklyanin:yz} and is called reflection algebra.
In the same way, $t(k)$ generates a subalgebra $\ct_{S}\subset 
\ck_{S}$, called in
analogy transmission algebra.

\section{From ${\overline\ca}_{S}$ to $\cc_{S}$}

Let $R(k)$ and $T(k)$ be any two matrix functions satisfying
\begin{equation}
R^\dagger(k)\, =\, R(-k)\, , \qquad T^\dagger(k) = T(k)
\end{equation}
\begin{equation}
T(k)\, T(k) + R(k)\, R(-k) = \II
\end{equation}
and the boundary Yang-Baxter equation
\begin{equation}
S_{12}(k_{1},k_{2})\, R_{1}(k_1)\, S_{21}(k_{2},-k_{1})\, R_{2}(k_2)\, =\,
R_{2}(k_2)\, S_{12}(k_{1},-k_{2})\, R_{1}(k_1)\,
S_{21}(-k_{2},-k_{1})
\end{equation}
It has been shown in \cite{Mintchev:2003ue} that $R(k)$ and $T(k)$ 
obey in addition
the {\it transmission} and {\it reflection-transmission} Yang-Baxter equations
\begin{eqnarray}
S_{12}(k_{1},k_{2})\, T_{1}(k_1)\, S_{21}(k_{2},k_{1})\, T_{2}(k_2)\, =\,
T_{2}(k_2)\, S_{12}(k_{1},k_{2})\, T_{1}(k_1)\, S_{21}(k_{2},k_{1}) \qquad \\
S_{12}(k_{1},k_{2})\, T_{1}(k_1)\, S_{21}(k_{2},k_{1})\, R_{2}(k_2)\, =\,
R_{2}(k_2)\, S_{12}(k_{1},-k_{2})\, T_{1}(k_1)\, S_{21}(-k_{2},k_{1})
\end{eqnarray}
together with
\begin{equation}
T(k)R(k)+R(k)T(-k) =  0
\end{equation}
This remarkable property of the pair $\{R(k),\, T(k)\}$
was discovered in \cite{Mintchev:2003ue} and is fundamental
in what follows. In fact,
by a direct calculation, starting from the EZF algebra ${\overline\ca}_{S}$,
one proves:
\begin{prop}
Let
\begin{eqnarray}
t(k) &=& L(k)T(k)L^{-1}(k) \label{t-L}\\
r(k) &=& L(k)R(k)L^{-1}(-k) \label{r-L}
\end{eqnarray}
where $L(k)$ is the well-bred operator of ${\overline\ca}_{S}$ and 
$R(k)$ and $T(k)$
are defined above. Then, $a(k)$, $a^\dag(k)$, $t(k)$ and
$r(k)$ obey the relations {\rm (\ref{rt-4})--(\ref{rt-12b})}.
\end{prop}
\begin{prop}
The map
\begin{equation}
\rho\left\{\begin{array}{lll}
a(k) & \rightarrow & \alpha(k)\, =\, t(k)a(k)+r(k)a(-k)\\
a^\dag(k) & \rightarrow &\alpha^\dag(k)
\, =\, a^\dag(k)t(k)+a^\dag(-k)r(-k)
\end{array}\right.
\end{equation}
extends to a homomorphism on ${\overline\ca}_{S}$.
\end{prop}
\prf
Direct calculation using the commutation relations of the ZF
algebra and the equations obeyed by $R(k)$ and $T(k)$.
\finprf
Let us remark that we have the identities
\begin{eqnarray}
a(k)&=& t(k)\alpha(k)+r(k)\alpha(-k)\\
a^\dag(k) &=& \alpha^\dag(k)t(k)+\alpha^\dag(-k)r(-k)
\end{eqnarray}
{\it However}, since $t(k)$ and $r(k)$ are still expressed in terms of $a(k)$
and $a^\dag(k)$, these are not really ``inversion formulas".

Note that $\alpha(k)$, $\alpha^\dag(k)$, $t(k)$ and $r(k)$ also obey
the relations (\ref{rt-4})--(\ref{rt-12b}).

The homomorphism $\rho$ is essential for constructing an RT algebra from
EZF. Indeed, we have:
\begin{theo}\label{theo1}
Let
\begin{eqnarray}
A(k) &=& \half\Big(a(k)+\alpha(k)\Big) =\half\Big(
[1+t(k)]a(k)+r(k)a(-k)\Big)\\
A^\dag(k) &=& \half\Big(a^\dag(k)+\alpha^\dag(k)\Big)=\half\Big(
a^\dag(k)[1+t(k)]+a^\dag(-k)r(-k)\Big)
\end{eqnarray}
Then $A(k)$, $A^\dag(k)$, $t(k)$ and $r(k)$ form a RT algebra.
\end{theo}
\prf
Direct calculation
\finprf
Note, {\it en passant},
that if one imposes on $R(k)$ the additional relation
$R(k)R(-k)=\II$, one gets $T(k)=0$ and one recovers the construction
\cite{Ragoucy:2001cf} of the boundary algebra $\cb_S$.
Let us also stress that $\rho$ is \underline{not} identity on the
algebra generated by $A(k)$, $A^\dag(k)$, $t(k)$ and
$r(k)$. However, we have:
\begin{theo}\label{theo34}
     The coset of ${\overline\ca}_{S}$ by the relation $\rho=id$ is a 
RT algebra.
\end{theo}
\prf
By coset by the relation $\rho=id$, we mean the coset by the two-sided
ideal generated by $Im(\rho-id)$: this coset is obviously an algebra. 
$A(k)$ and
$A^\dag(k)$ defined in theorem \ref{theo1} can be taken as
representatives of a generating system for the coset which is thus
homomorphic to the RT algebra.
\finprf
To conclude this section, let us remark that, in the EZF algebra, one
has the identities:
\begin{eqnarray}
A(k)&=& t(k)A(k)+r(k)A(-k)\\
A^\dag(k) &=& A^\dag(k)t(k)+A^\dag(-k)r(-k)
\end{eqnarray}
These identities show that 
for the RT algebra constructed in theorem \ref{theo1}, the 
reflection-transmission automorphism (as it was introduced in \cite{Mintchev:2003ue}) 
is indeed the identity, in agreement with theorem \ref{theo34}. 
However, this automorphism (defined only on 
the RT algebra) has not to be 
confused with the above homomorphism $\rho$, defined on the whole EZF 
algebra.

\section{Hamiltonians and their symmetries}

It is known that one can associate to any ZF algebra a natural hierarchy
of Hamiltonians
\begin{equation}
H^{(n)}_{ZF}=\int_{\RR}dk\, k^n\, a^\dag(k)a(k),\ n\in\ZZ_{+}
\end{equation}
They obey to $[H^{(n)}_{ZF},H^{(m)}_{ZF}]=0$  (so that they can
indeed be considered as Hamiltonians) and
admit as symmetry algebra the quantum group generated by the
well-bred operators:
\begin{equation}
[H^{(n)}_{ZF},L(k)]=0,\ \forall n
\end{equation}
   Note that we have the identity
\begin{equation}
H^{(2n)}_{ZF}=\int_{\RR}dk\, k^{2n}\, a^\dag(k)a(k)
=\int_{\RR}dk\, k^{2n}\, \alpha^\dag(k)\alpha(k)
\label{eqH2n}
\end{equation}
which shows that $\rho(H^{(2n)}_{ZF})=H^{(2n)}_{ZF}$, so that 
$H^{(2n)}_{ZF}$ survives the coset of ${\overline\ca}_{S}$ by $\rho=id$.

In a similar way, for any RT algebra, we introduce
\begin{equation}
H_{RT}^{(n)}=\int_{\RR}dk\, k^n\, A^\dag(k)A(k)
\end{equation}
These Hamiltonians can be viewed as the representative of $H^{(n)}_{ZF}$
in the coset algebra of theorem \ref{theo34}. 
They obey to the following relations:
\begin{equation}
[H^{(m)}_{RT},H^{(n)}_{RT}]= [(-1)^m - (-1)^n] \int_{\RR}dk\, k^{m+n}\, \,
A^\dag(k)r(k)A(-k) , \label{comHRT}
\end{equation}
It shows that Hamiltonians of the same parity form a
commutative algebra, and thus define a hierarchy corresponding to the
integrable systems with impurity,
in accordance with eq. (\ref{eqH2n}). They
 act on $A(k)$ and $A^\dagger(k)$ as
\begin{eqnarray}
{[H_{RT}^{(n)},A(k)]} &=&
-k^n\Big( [\II+t(k)]A(k)+(-1)^nr(k)A(-k)\Big)\\
{[H_{RT}^{(n)},A^\dag(k)]} &=&
k^n\Big( A^\dag(k)[\II+t(k)]+(-1)^nA^\dag(-k)r(-k)\Big)
\end{eqnarray}
Then, by a direct calculation and {\it without using the embedding}
$\cc_{S}\subset {\overline\ca}_{S}$, one proves:
\begin{prop}
The subalgebra $\ck_{S}$ is a symmetry algebra of the hierarchy
$H_{RT}^{(n)}$:
\begin{equation}
[H^{(n)}_{RT},t(k)]=0 \mb{and} [H^{(n)}_{RT},r(k)]=0
\end{equation}
\end{prop}
This result provides an universal model-independent description of
the symmetry content of the hierarchy $H_{RT}^{(n)}$.
It is quite remarkable that $r(k)$ and $t(k)$ encode both the particle-impurity
interactions and the quantum integrals of motion. This property of 
$\ck_{S}$ has been already
applied \cite{Mintchev:2001aq} with success for studying the symmetries of
the $gl(N)$-invariant nonlinear Schr\"odinger equation on the half line.
 
{}Finally, using the embedding of the RT algebra into ${\overline 
\ca}_{S}$, one shows:
\begin{prop}
When considering the embedding $\cc_{S}\subset {\overline \ca}_{S}$,
the subalgebra $\ck_{S}$ becomes a Hopf coideal of the quantum group
$\cu_{S}$, \ie in algebraic terms one has $\ck_{S}\subset\cu_{S}$ and
\begin{equation}
\Delta(\ck_{S})\subset \cu_{S}\otimes \ck_{S}
\end{equation}
where $\Delta$ is the the coproduct of $\cu_{S}$.
\end{prop}
\prf The construction (\ref{t-L}) and
(\ref{r-L}) ensures the algebra embedding. Using the coproduct of
$\cu_{S}$, one gets
\begin{eqnarray}
\Delta[t_{ab}(k)] &=& \sum_{x,y} L_{ax}(k) L^{-1}_{yb}(k)\otimes
t_{xy}(k)\\
\Delta[r_{ab}(k)] &=& \sum_{x,y} L_{ax}(k) L^{-1}_{yb}(-k)\otimes
r_{xy}(k)
\end{eqnarray}
\finprf

Let us observe in conclusion that there is a simple relation between 
the hierarchies
with and without impurity, which reads
\begin{equation}
H_{RT}^{(n)}=H_{ZF}^{(n)}+ \int_{\RR}dk\,k^n\,
a^\dag(k)\Big(r(k)a(-k)+t(k)a(k)\Big)
\label{relation}
\end{equation}
Eq. (\ref{relation}) generalizes the result of \cite{ZFconf} for the 
boundary algebra
${\cal B}_S$.

\section{Conclusions}

The results of this paper clarify the connection between EZF and RT algebras.
We have shown above that a suitable coset algebra,
constructed in terms of the EZF algebra, is homomorphic to the RT 
algebra.
This feature provides a precise mathematical meaning of
the physical observation \cite{Mintchev:2002zd, Mintchev:2003ue}
that the introduction of an impurity,
preserving the integrability of a system, is equivalent to
imposing some supplementary constraints on the system. The latter
implement the consistency between the scattering in the bulk and
the interaction with the impurity.

\end{document}